\documentclass[11pt,a4paper]{article}
\usepackage{bbm}
\usepackage{mathrsfs}
\usepackage{amsfonts}
\usepackage{}

\usepackage{graphicx}
\usepackage{amsfonts,amsmath, amssymb}

\newtheorem{theo}{\textbf{\ \ \quad Theorem}}[section]

\newtheorem{lem}{\textbf{\ \ \quad Lemma}}[section]
\newtheorem{remark}{\textbf{\ \ \quad Remark}}[section]

\newtheorem{prop}{\textbf{\ \ \quad Proposition}}[section]

\newcommand{\lbl}[1]{\label{#1}}

\newcommand{\be}{\begin{equation}}
\newcommand{\ee}{\end{equation}}
\newcommand\bes{\begin{eqnarray}}
\newcommand\ees{\end{eqnarray}}
\newcommand{\bess}{\begin{eqnarray*}}
\newcommand{\eess}{\end{eqnarray*}}

\newcommand{\nm}{\nonumber}

\newcommand{\ds}{\displaystyle}

{\setlength\arraycolsep{2pt}}

\title{The effect of noise intensity on stochastic parabolic equations}

\author{Guangying Lv\footnote{Institute of Applied Mathematics, Henan University, Kaifeng, Henan 475001, China;    
and Institute of Mathematics, School of Mathematical Science, Nanjing Normal University, Nanjing 210023, China.  
E-mail: gylvmaths@henu.edu.cn. Research supported in part by NSFC of China grants 11771123, 11171064.
} \quad Hongjun Gao\footnote{Institute of Mathematics, School of Mathematical Science, 
Nanjing Normal University, Nanjing 210023, China. E-mail: gaohj@njnu.edu.cn} \quad  Jinlong Wei\footnote{
School of Statistics and Mathematics, Zhongnan University of Economics and Law, Wuhan, Hubei 430073, China. 
E-mail:  weijinlong.hust@gmail.com} \quad Jiang-Lun Wu\footnote{Department of Mathematics, Swansea University, 
Swansea SA2 8PP, UK.  E-mail: j.l.wu@swansea.ac.uk \,\, http://www.swansea.ac.uk/staff/science/maths/j.l.wu/}}

\begin{document}
\maketitle

\medskip

\begin{abstract} In the present paper, the effect of noise intensity on stochastic 
parabolic equations is discussed. We focus on the effect of noise on the energy 
solutions of the stochastic parabolic equations. By utilising It\^o's formula and 
the energy estimate method, we obtain the excitation indices of the energy solutions  
$u$ at any finite time $t$. Furthermore, we improve certain existing results in the 
literature by presenting a comparably simple method to show those existing results.



\end{abstract}

\baselineskip=15pt

\section{Introduction}
\setcounter{equation}{0}
In recent years, many authors attempt to explore the role of the noise in various dynamical equations in both analytical and numerical aspects.
For example, noise can make the solution smooth \cite{FGP}, can prevent singularities in linear transport equations \cite{FF}, can prevent collapse
of Vlasov-Poisson point charges \cite{DFV}, and also can induce singularities (finite time blow up of solutions) \cite{C2009,C2011,LD}. In the
present paper, we focus on the effect of noise on parabolic equations.

The concept of  ``Intermittency" is the property that the solution $u(t,x)$ develops extreme oscillations at certain values of $x$, typically when $t$
is going to be large. Intermittency was announced first (1949) by Batchelor and Townsend in a WHO conference in Vienna \cite{BT}, and slightly
later by Emmons \cite{E} in the context of boundary layer turbulence. Meanwhile, intermittency has been observed in an enormous number of
scientific disciplines. For example, intermittency is observed as ``spikes" and ``shocks" in neural activity and in finance,
respectively. Tuckwell \cite{T} contains a gentle introduction to SPDEs in neuroscience.

Recently, Khoshnevisan-Kim in \cite{KK1,KK2} considered the following stochastic heat equation
\bes
\frac{\partial}{\partial t}u(t,x)=\mathcal {L}u(t,x)+\lambda\sigma(u(t,x))\xi(t,x) 
\lbl{1.1}\ees
where $t\in(0,\infty)$ stands for the time variable, $x\in G$ the space variable with $G$ being  a given nice state space, such as 
$\mathbb{R},\,\mathbb{Z}$ (a discrete set) or a finite interval like $[0,1]$, and the initial data value $u_0:\,G\rightarrow\mathbb{R}$
is deterministic (i.e., non random) and is well behaved. The operator $\mathcal {L}$ acts on the spatial variable $x\in G$ only, and 
is taken to be the generator of a nice Markov stochastic process on $G$, and $\xi$ denotes \text{\emph{space-time white noise}} 
on $(0,\infty)\times G$. Whereby, $\lambda>0$ is a constant and the coefficient 
$\sigma:\mathbb{R}\to\mathbb{R}$ is supposed to be a Lipschitz continuous function.

Let $u$ be a mild solution of (\ref{1.1}) with given initial data $u_0$. Set 
$$u(t):=u(t,\cdot):D\to\mathbb{R},\quad t\in[0,\infty)$$  
and then define
 \bes
\mathscr{E}_t(\lambda):=\sqrt{\mathbb{E}\left(\|u(t)\|_{L^2(G)}^2\right)},
\lbl{1.2}\ees
which stands for the energy of the solution at time $t$. In papers \cite{KK1,KK2,MM,MTL}, 
the authors showed that the energy $\mathscr{E}_t(\lambda)$ behaves like 
$exp({\rm const}\cdot\lambda^q)$, for certain fixed positive constant $q$, as $\lambda\uparrow\infty$. In order to do so, 
the following two quantities have been introduced  
   \bes
   \underline{\mathbbm{e}}(t):=\lim\limits\inf_{\lambda\uparrow\infty}\frac{\log\log\mathscr{E}_t(\lambda)}{\log\lambda},\ \
   \bar{\mathbbm{e}}(t):=\lim\limits\sup_{\lambda\uparrow\infty}\frac{\log\log\mathscr{E}_t(\lambda)}{\log\lambda}.
    \lbl{1.3} \ees
Clearly, $\underline{\mathbbm{e}}$ and $\bar{\mathbbm{e}}$ represent the lower and upper excitation indices of $u$ at
time $t$, respectively. In many interesting cases, $\underline{\mathbbm{e}}(t)$ and $\bar{\mathbbm{e}}(t)$ are exactly equal, 
and they do not depend on the time variable $t\in[0,\infty)$. In such situations, we tacitly write $\mathbbm{e}$
for that common value, just for simplicity.

In paper \cite{KK1}, Khoshnevisan-Kim proved that

(i) If $G$ is discrete, then $\bar{\mathbbm{e}}(t)\leq2$ for
all $t\geq0$. Furthermore, it hold that $\bar{\mathbbm{e}}=2$ if
   \bes
l_\sigma:=\inf_{z\in\mathbb{R}\backslash\{0\}}|\frac{\sigma(z)}{z}|>0.
   \lbl{1.4}\ees

(ii) Suppose that $G$ is connected and (\ref{1.4}) holds, then $\underline{\mathbbm{e}}(t)\geq4$ for
all $t\geq0$, provided that in addition either $G$ is non compact or $G$ is compact, metrizable, and
has more than one element.

(iii) For every $\theta\geq4$ there exist models of the triple $(G,\mathcal {L},u_0)$ for
which $\mathbbm{e}=\theta$. One such model is that $\mathcal {L}:=-(-\Delta)^{\frac{\alpha}{2}}$
(the generator of a symmetric $\alpha$-stable L\'{e}vy process) for $1<\alpha\leq2$.

In \cite{KK2}, Khoshnevisan-Kim considered the following problem for the stochastic evolution equation  
 \bes\left\{\begin{array}{llllll}
 \ds\frac{\partial}{\partial t}u(t,x)=\frac{\partial^2}{\partial x^2}u(t,x)+\lambda\sigma(u(t,x))\dot{w}(t,x), \ \ \
 &0<x<L,\,t>0,\\
u(t,0)=u(t,L)=0, \ \ &t>0,\\
u(0,x)=u_0(x),
 \end{array}\right.\lbl{1.5}\ees
where $\dot{w}$ is a space-time white noise, $L>0$ is fixed, $u_0(x)\geq0$ is a non-random,  bounded continuous function and 
$\sigma:\mathbb{R}\rightarrow\mathbb{R}$ is a Lipschitz continuous function with $\sigma(0)=0$. Let
   \bes
l_\sigma:=\inf_{z\in\mathbb{R}\backslash\{0\}}|\frac{\sigma(z)}{z}|>0,\ \
L_\sigma:=\sup_{z\in\mathbb{R}\backslash\{0\}}|\frac{\sigma(z)}{z}|>0.
   \lbl{1.6}\ees
They derived the following
   \bess
\frac{l_\sigma^2t}{2}\leq\liminf_{\lambda\rightarrow\infty}\frac{1}{\lambda^2}\log\mathscr{E}_t(\lambda),\ \
\limsup_{\lambda\rightarrow\infty}\frac{1}{\lambda^4}\log\mathscr{E}_t(\lambda)\leq 8L^4_\sigma t.
  \eess
More recently, Foondun-Joseph \cite{MM} complemented the results of \cite{KK2}, that is, they obtained
$\mathbbm{e}=4$. It is easy to see that a mild solution $u$ of (\ref{1.5}) which is adapted to the natural filtration of the white noise $\dot{w}$ 
and satisfies the following mild formulation of the evolution equation
   \bes
u(t,x)=(\mathcal {G}_Du)(t,x)+\lambda\int_0^t\int_0^Lp_D(t-s,x,y)\sigma(u(s,y))w(dsdy),
  \lbl{1.7}\ees
where
  \bess
(\mathcal {G}_Du)(t,x):=\int_0^Lu_0(y)p_D(t,x,y)dy,
   \eess
and $p_D(t,x,y)$ denotes the Dirichlet heat kernel, $D:=[0,L]$. They used the estimate of
kernel $p_D(t,x,y)$ and a new Gronwall's inequality to prove that $\mathbbm{e}=4$. Using similar
method, Liu-Tian-Foondun \cite{MTL} considered the fractional Laplacian on a bounded domain.

Now, given a complete probability space endowed with a filtration $(\Omega,\mathcal {F},\{\mathcal {F}_t\}_{t\geq0},\mathbb{P})$, let us consider 
the following linear stochastic differential equation (SDE) with $\lambda>0$ being given before 
   \bess
dX_t=\lambda X_tdB_t,\ \ t>0, \ \ X_0=x\in D.
  \eess
For simplicity, we assume that $B(t)$ is a standard one-dimensional Brownian motion on $(\Omega,\mathcal {F},\{\mathcal {F}_t\}_{t\geq0},\mathbb{P})$.
It is easy to see that the unique solution of the above SDE is explicitly given by
   \bess
X_t=xe^{-\frac{\lambda^2}{2}t}e^{\lambda B(t)}, \ \ \ \ t\in[0,+\infty).
   \eess
Direct calculations then show that
\bess
&&\mathbb{E}[X_t]=xe^{-\frac{\lambda^2}{2}t}e^{\frac{\lambda^2}{2}t}=x;\quad \mathbb{E}[X^2_t]=x^2e^{-\lambda^2t}e^{2\lambda^2t}=x^2e^{\lambda^2t};\\
&&\mathbb{E}[X^p_t]=x^pe^{-\frac{\lambda^2p}{2}t}e^{\frac{\lambda^2p^2}{2}t}=
x^pe^{\frac{\lambda^2p(p-1)}{2}t}
\eess
for $p>1$. This then implies that for $p>1$
   \bess
\lim\limits_{\lambda\to\infty}\frac{\log\log\left(\mathbb{E}[X^p_t]\right)}{\log\lambda}=2,
   \eess
which yields that the excitation index of $X_t$ is $2$. This is clearly
different from the results obtained in \cite{MM,KK2,Xie2016},
where the authors proved the excitation indice of $u(t,x)$ of (\ref{1.5}) is $4$
for $x\in[\epsilon,L-\epsilon]$ ($\epsilon$ is a sufficiently small constant). A natural 
and very interesting question then appeared to be that is there some kind of solutions of (\ref{1.5}) with the associated indices being $2$? 
This motivates us to initiate the present paper. 

Another propose of our paper is to introduce a comparably simpler method to prove the result of \cite{MM} in a simple case. That is, we 
consider the following stochastic parabolic equations
    \bes\left\{\begin{array}{llllll}
du(t,x)=\Delta u(t,x)dt+\lambda u(t,x)dB_t, \ \ \
 &x\in D,\,t>0,\\[1mm]
u|_{\partial D}=0, \ \ &t>0,\\
u(0,x)=u_0(x),
 \end{array}\right.\lbl{1.8}\ees
where $D\subset\mathbb{R}^n$ ($n\geq1$), $B_t$ is a standard one-dimensional Brownian
motion on $(\Omega,\mathcal {F},\{\mathcal {F}_t\}_{t\geq0},\mathbb{P})$ as given above.
We will obtain the similar result to \cite{MM} by changing the stochastic parabolic
equations into random parabolic equations. Moreover,
it is not hard to find that, in the earlier results, the authors only consider the expression 
$\sqrt{\mathbb{E}\left(\|u(t)\|_{L^2(G)}^2\right)}$, and we can consider the expression 
$\left[\mathbb{E}\left(\|u(t)\|_{L^p(G)}^p\right)\right]^{1/p}$, $p>0$, which is clearly an interesting generalisation.

In this paper, we will focus on the noise excitability of energy solution for
some parabolic equations. We obtain a new result regarding the noise excitability, that is,
$\mathbbm{e}=2$ under the same condition as in \cite{MM} when the noise is only the time 
white noise  (not the space-time white noise). The contribution of our paper is that we 
consider energy solutions (comparing to \cite{MM} where mild solutions are considered). 

The rest of the paper is organised as follows. In Section 2, some preliminaries and main results are given. Section 3 is devoted to the 
proofs of the main results. In Section 4, we consider a special noise case of (\ref{1.7}) and we discuss noise excitability of stochastic 
equations involving nonlocal operators. 

\section{Preliminaries and two main results}
\setcounter{equation}{0}

Inspired by \cite{KK1,KK2,MM,eP,Ta}, in this paper, we consider the simple case
 \bes\left\{\begin{array}{llllll}
 du(t,x)=\Delta u(t,x)dt+\lambda\sigma(u(t,x))dB_t, \ \ \
 &x\in D,\,t>0,\\[1mm]
u|_{\partial D}=0, \ \ &t>0,\\
u(0,x)=u_0(x),\ & x\in D
 \end{array}\right.\lbl{2.6}\ees
where $D\subset\mathbb{R}^n$ ($n\geq1$) is a bounded domain, and $B_t$ denotes one dimensional Brownian motion.
The existence of solutions of (\ref{2.6}) was obtained by \cite{Cb}.

The main results of this paper are formulated as the following
\begin{theo}\lbl{t2.1} Assume that (\ref{1.6}) holds.
The noise excitation index of the energy solution to (\ref{2.6}) with initial data $u_0(x)\geq0 (\not\equiv0)$ is $2$.
  \end{theo}

If the one-dimensional Brownian motion is replaced by $Q$-Wiener process, where $Q$ is a trace 
class operator on $L^2(D)$, the result of Theorem \ref{t2.1} still holds. In fact, for the following equation 
(i.e., the case driven by time white noise), we have the following result.
 \bes\left\{\begin{array}{llllll}
 \frac{\partial}{\partial t}u(t,x)=\Delta u(t,x)+\lambda\sigma(u(t,x))\dot w(t,x) \ \ \
 &x\in D,\,t>0,\\[1mm]
u|_{\partial D}=0, \ \ &t>0,\\
u(0,x)=u_0(x),\ & x\in D
 \end{array}\right.\lbl{2.7}\ees
where $D\subset\mathbb{R}^n$ ($n\geq1$) is a bounded domain and $w(t,x)$ is a $Q$-Wiener process.
The existence of solutions of (\ref{2.7}) was also obtained by \cite{Cb}.

\begin{theo}{\lbl{t2.2}} Assume that (\ref{1.6}) holds. Let $w(t,x)$ be a $Q$-Wiener process with covariance 
    \bess
\mathbb{E}[w(t,x)w(s,y)]=(t\wedge s)q(x,y),\quad s,t\in(0,+\infty), \, x,y\in D 
\eess
where $q:D\times D\to\mathbb{R}$ is the kernel of the trace class operator $Q:L^2(D)\to L^2(D)$. 
   
In this case, the noise $\dot{w}(t,x)$ is white in time and colored in the space variable.   
Assume that $0<\sup_{x\in D}q(x,x)\leq q_1<\infty$, then, the upper excitation index of the solution to (\ref{2.7}) with 
initial data $u_0(x)\geq0$ is $2$. Furthermore, if $\sigma\geq 0$ (or $\leq0$) and there is a positive real number $q_0>0$ 
such that $q_0<\inf_{x,y\in D}q(x,y)$, then the excitation index of the solution to (\ref{2.7}) with initial data 
$u_0(x)\geq0,\not\equiv0$ is $2$.
\end{theo}

\begin{remark}\lbl{r2.1} Now, we give the reason why we can not consider
the case that $g(u)$ satisfies \textbf{\emph{local Lipschitz}} condition. More precisely, consider the following general case
 \bes\left\{\begin{array}{llll}
du(t)=[Au(t)+f(u(t))]dt+\lambda \sigma(u(t))dw(t),\ \ \ \ t>0,\\
u(0)=u_0,
   \end{array}\right.\lbl{2.8}\ees
where $A$ is a divergence operator, $f$ and $\sigma$ satisfy the local Lipschitz condition.
For example, let $f(u)\geq au^{1+\alpha}$ and $\sigma(u)=u^m$. Then the solutions of {\rm(\ref{2.8})} will blow up in finite
time (see \cite{C2011,LD}). Moreover, the largest existence time $T\rightarrow0$ as $\lambda\rightarrow\infty$. So
we cannot consider problem {\rm(\ref{2.8})}.
\end{remark}

\section{The proofs of our main results}
\setcounter{equation}{0}
In this section, we will prove Theorem \ref{t2.1} and Theorem \ref{t2.2} by using energy method.
Let us first prove Theorem \ref{t2.1}.

{\bf Proof of Theorem \ref{t2.1}.}  By using the idea of \cite{eP,Ta}, one can prove that
there exists a unique energy solution. It follows from the results of \cite{LD} that the energy solution
will keep positive if the initial data $u_0\geq0$ almost surely. We divide our proof into two steps.

Step 1: $\bar{\mathbbm{e}}(t)=2$.

By It\^{o} formula, we have
   \bes
\|u(t)\|_{L^2}^2&=&\|u_0\|_{L^2}^2+2\int_0^t\langle \Delta u(s,x),u(s,x)\rangle
ds+2\lambda\int_0^t\int_Du(s,x)\sigma(u(s,x))dxdB_s\nm\\
&&
+\lambda^2\int_0^t\int_D\sigma^2(u(s,x))dxds.
   \lbl{3.1}\ees
Integrating by parts shows that
     \bess
\|u(t)\|_{L^2}^2&=&\|u_0\|_{L^2}^2-2\int_0^t\|\nabla u(s)\|_{L^2}^2ds+2\lambda\int_0^t\int_Du(s,x)\sigma(u(s,x))dxdB_s\nm\\
&&
+\lambda^2\int_0^t\int_D\sigma^2(u(s,x))dxds\\
&\leq&\|u_0\|_{L^2}^2+2\lambda\int_0^t\int_Du(s,x)\sigma(u(s,x))dxdB_s
+\lambda^2\int_0^t\int_D\sigma^2(u(s,x))dxds,
    \eess
which implies
      \bess
\mathbb{E}\|u(t)\|_{L^2}^2&\leq&\mathbb{E}\|u_0\|_{L^2}^2
+\lambda^2\mathbb{E}\int_0^t\int_D\sigma^2(u(s,x))dxds\\
&\leq&\mathbb{E}\|u_0\|_{L^2}^2
+L_\sigma\lambda^2\int_0^t\mathbb{E}\|u(s)\|_{L^2}^2dxds.
    \eess
It follows from Gronwall's inequality that
   \bes
\mathbb{E}\|u(t)\|_{L^2}^2\leq \mathbb{E}\|u_0\|_{L^2}^2e^{L_\sigma\lambda^2t},
   \lbl{3.2}\ees
which implies that $\bar{\mathbbm{e}}(t)\leq2$.

Step 2: $\underline{\mathbbm{e}}(t)=2$.

In order to get the lower bounded, let us consider the following eigenvalue problem for the elliptic equation
  \bes\left\{\begin{array}{llll}
-\Delta \phi=\lambda \phi, \ \ \ \ \ \ \  \ {\rm in} \  D,\\
\phi=0, \ \ \qquad\ \qquad  {\rm on}\ \partial D.
   \end{array}\right.\lbl{3.3}\ees
Since all the eigenvalues are strictly positive, increasing and
the eigenfunction $\phi$ corresponding to the smallest eigenvalue
$\lambda_1$ does not change sign in domain $ D$, as shown in \cite{GT}, 
then one can normalise it in such a way that
   \bess
\phi(x)>0 \ \mbox{in} \ D, \ \ \int_ D \phi(x)dx=1.
   \eess

Noting that under the assumptions of Theorem \ref{t2.1}, the solutions
of (\ref{1.5}) will remain positive, thus we can consider $(u(t),\phi)$ due to 
the fact that $(u(t),\phi)>0$. Denote $\hat u(t):=(u(t),\phi)$.
By applying It\^{o}'s formula to $\hat u^2(t)$ and making use of (\ref{3.3}), we get
   \bes
\hat u^2(t)&=&(u_0,\phi)^2-2\lambda_1\int_0^t\hat u^2(s)ds
+2\lambda\int_0^t\int_ D \hat u_(s)\sigma(u_s(x))\phi(x)dxdB_s\nm\\
&&+\lambda^2\int_0^t\int_ D\sigma^2(u_s(x))\phi^2(x)dxds\nm\\
&\geq&(u_0,\phi)^2-2\lambda_1\int_0^t\hat u^2(s)ds
+2\lambda\int_0^t\int_ D \hat u(s)\sigma(u(s,x))\phi(x)dxdB_s\nm\\
&&+\lambda^2l^2_\sigma\int_0^t\int_ Du^2(s,x)\phi^2(x)dxds\nm\\
&\geq&(u_0,\phi)^2-2\lambda_1\int_0^t\hat u^2(s)ds
+2\lambda\int_0^t\int_ D \hat u_(s)\sigma(u(s,x))\phi(x)dxdB_s\nm\\
&&+\lambda^2l^2_\sigma\int_0^t\hat u^2(s)ds.
   \lbl{3.4}\ees
Taking mean norm then yields that
   \bes
\mathbb{E}\hat u^2(t)
&\geq&\mathbb{E}(u_0,\phi)^2-2\lambda_1\int_0^t\mathbb{E}\hat u^2(s)ds
+\lambda^2l^2_\sigma\int_0^t\mathbb{E}\hat u^2(s)ds.
   \lbl{3.5}\ees
By the comparison principle, we know that
   \bess
\mathbb{E}\hat u^2(t)\geq \mathbb{E}(u_0,\phi)^2e^{(\lambda^2l^2_\sigma-2\lambda_1)t}.
   \eess
Due to
  \bess
\hat u^2(t)=(u,\phi)^2\leq \|\phi\|_{L^\infty}\|u\|_{L^2}^2,
   \eess
we have $\underline{\mathbbm{e}}(t)\geq2$. So we
have $\mathbbm{e}=2$. $\Box$

{\bf Outline of the proof of Theorem \ref{t2.2}.}  Similar to the proof of Theorem \ref{t2.1},
equation (\ref{2.6}) has a unique positive energy solution.

From (\ref{3.1}), we have
  \bess
\|u(t)\|_{L^2}^2&=&\|u_0\|_{L^2}^2+2\int_0^t\langle \Delta u(s,x),u(s,x)\rangle
ds+2\lambda\int_0^t\int_Du(s,x)\sigma(u(s,x))w(dxds)\nm\\
&&+
\lambda^2\int_0^t\int_Dq(x,x)\sigma^2(u(s,x))dxds\\
&\leq&\|u_0\|_{L^2}^2+2\lambda\int_0^t\int_Du(s,x)\sigma(u(s,x))w(dxds)
+q_1\lambda^2L_\sigma\int_0^t\int_Du^2(s,x)dxds.
    \eess
Then taking expectation on both sides and using Gr\"{o}nwall's inequality, we have $\bar{\mathbbm{e}}(t)\leq2$.

Similar to the proof of Theorem \ref{t2.1}, we have further
\bes
\hat u^2(t)&=&(u_0,\phi)^2-2\lambda_1\int_0^t\hat u^2(s)ds+2\lambda\int_0^t\int_ D \hat u(s)\sigma(u(s,x))\phi(x)w(dx,ds)\nm\\
&&+\lambda^2\int_0^t\int_ D\int_ D \sigma(u(s,x))\phi(x)q(x,y)\sigma(u(s,y))\phi(y)dxdyds
\nm\\
&\geq&(u_0,\phi)^2-2\lambda_1\int_0^t\hat u^2(s)ds+2\lambda\int_0^t\int_ D \hat u(s)\sigma(u(s,x))\phi(x)w(dx,ds)\nm\\
&&+\lambda^2q_0\int_0^t\int_ D\int_ D  \sigma(u(s,x))\phi(x)\sigma(u(s,y))\phi(y)dxdyds\nm\\
&\geq&(u_0,\phi)^2-2\lambda_1\int_0^t\hat u^2(s)ds+2\lambda\int_0^t\int_ D \hat u(s)\sigma(u(s,x))\phi(x)w(dx,ds)\nm\\
&&+\lambda^2q_0l^2_\sigma\int_0^t\hat u^2(s)ds,
\lbl{3.6}\ees
which implies that $\underline{\mathbbm{e}}(t)\geq2$. So we
have $\mathbbm{e}=2$. $\Box$

\begin{remark}\lbl{r2.1}1. We have considered the problem with higher space dimensions 
as we study the equations perturbed by a noise white in time and colored in the space 
variable. While in papers \cite{KK1,KK2,MM,MTL},
the authors only considered one space dimension due to their equations are driven by 
space time white noise.

2. The Laplace operator $\Delta$ can be substituted by the divergent operator $A$. 
\end{remark}

\section{A special case and the noise excitability for nonlocal operators}
\setcounter{equation}{0}
In this section, we consider the following problem 
    \bes\left\{\begin{array}{llllll}
du(t,x)=\Delta u(t,x)dt+\lambda u(t,x)dB_t, \ \ \
 &x\in D,\,t>0,\\[1mm]
u|_{\partial D}=0, \ \ &t>0,\\
u(0,x)=u_0(x),\ \ &x\in D,
 \end{array}\right.\lbl{4.1}\ees
where $D\subset\mathbb{R}^n$ ($n\geq1$), $B_t$ is a standard one-dimensional Brownian
motion on a stochastic basis $(\Omega,\mathcal {F},\{\mathcal {F}_t\}_{t\geq0},\mathbb{P})$.

We first give a equivalent equation to (\ref{4.1}).
\begin{lem}\lbl{l4.1}
Let $u$ be a weak solution of (\ref{4.1}). Then the function $v$ defined via 
   \bess
v(t,x):=e^{-\lambda B_t}u(t,x), \ \ t>0,\ x\in D
   \eess
solves the following deterministic equations for random filed $v(t,x)$  
\bes\left\{\begin{array}{llllll}
\frac{\partial}{\partial t}v(t,x)=\Delta v(t,x)-\frac{\lambda^2}{2} v(t,x), \ \ \
&x\in D,\,t>0,\\[1mm]
v|_{\partial D}=0, \ \ &t>0,\\
v(0,x)=u_0(x),\ \ &x\in D.
\end{array}\right.\lbl{4.2}\ees
\end{lem}
The proof of this lemma is standard, see e.g. the proof of Proposition
1.1 of \cite{DL}. We therefore omit it here.

\begin{theo}\lbl{t4.1} Let $u$ be a weak solution of (\ref{4.1}) with deterministic 
initial data $u_0$ satisfying
    \bes
c_1\leq u_0(x)\leq c_2,\ \ \ \ \forall x\in D,
  \lbl{4.3}\ees
where $c_i$, $i=1,2$, are positive constants. Then we have, for $p>0$,
   \bes
2\leq\liminf_{\lambda\to\infty}\frac{\log\log\mathcal {E}_t(\lambda)}{\log\lambda}\leq \limsup_{\lambda\to\infty}\frac{\log\log\mathcal {E}_t(\lambda)}{\log\lambda}
\leq 2,
   \lbl{4.4}\ees
where $\mathcal {E}_t(\lambda):=\left[\mathbb{E}\left(\|u(t)\|_{L^p(D)}^p\right)\right]^{1/p}$.
 \end{theo}

{\bf Proof.} It follows from Lemma \ref{l4.1} that the solutions of (\ref{4.1}) can be expressed  as 
   \bess
u(t,x)=e^{\lambda B_t}v(t,x).
   \eess
It follows from the classical parabolic theory that the solutions $v$ of (\ref{4.2}) can be written as
   \bess
v(t,x)&=&e^{\frac{\lambda^2}{2}t}(e^{t\Delta}u_0)(x)\\
&=&e^{\frac{\lambda^2}{2}t}\int_Dp_D(t,x-y)u_0(y)dy, \ \ \ a.s.,
   \eess
where $p_D(t,x)$ is the kernel function of the Dirichlet Laplacian $\Delta$ on $D$. By using
(\ref{4.3}), we have
    \bess
\hat c_1e^{\frac{\lambda^2}{2}t}\leq v(t,x)\leq \hat c_2e^{\frac{\lambda^2}{2}t}, \ \ \ a.s.,
   \eess
which implies that
    \bess
\mathbb{E}\left[\|u(t)\|_{L^p(D)}^p\right]
&=&\mathbb{E}\left[\|v(t)e^{\lambda B_t}\|_{L^p(D)}^p\right]\\
&\geq& \tilde c_1e^{\frac{\lambda^2}{2}pt}\mathbb{E}\left[e^{\lambda p B_t}\right]\\
&=&\tilde c_1e^{\frac{\lambda^2}{2}pt}e^{\lambda^2p^2t}
   \eess
and
    \bess
\mathbb{E}\left[\|u(t)\|_{L^p(D)}^p\right]
&=&\mathbb{E}\left[\|v(t)e^{\lambda B_t}\|_{L^p(D)}^p\right]\\
&\leq& \tilde c_2e^{\frac{\lambda^2}{2}pt}\mathbb{E}\left[e^{\lambda p B_t}\right]\\
&=&\tilde c_2e^{\frac{\lambda^2}{2}pt}e^{\lambda^2p^2t}.
   \eess
Combining the above two inequalities, we get the desired result. The proof is thus complete. $\Box$

Next, we will consider the following initial value problem for nonlocal equations 
\bes\left\{\begin{array}{lll}
\frac{\partial}{\partial t}u(t,x)=-(-\Delta)^{\frac{\alpha}{2}}u(t,x)+\lambda\sigma(u(t,x))\dot{w}(t,x), \ \ &x\in \mathbb{R},\ t>0\\[1.5mm]
u(0,x)=u_0(x), \ \ \ &x\in \mathbb{R},
\end{array}\right.\lbl{4.6}\ees
where $\alpha\in(1,2]$, $(-\Delta)^{\frac{\alpha}{2}}$ is the $L^2$-generator of a symmetric
$\alpha$-stable process $X_t$ such that $\mathbb{E}\exp(i\xi\cdot X_t)=\exp(-t|\xi|^\alpha)$,
$\{\dot{w}(x,t)\}_{t\geq0,x\in\mathbb{R}}$
denotes the space-time white noise. In paper \cite{MTL}, the authors considered
the equation (\ref{4.6}) on bounded domain. Here we would like to generalise
the result to the situation of whole spatial space. 

When $\sigma$ satisfies global Lipschitz continuous condition, it is routine to show that (\ref{4.6}) has a unique global mild solution, 
see e.g. the monographs \cite{Cb,PZ,LiuRockner}, as well as Dalang \cite{R} and Foondun-Khoshnevisan \cite{FK1}. 
It is easy to see that the mild solution of (\ref{4.6}) fulfills the following mild formulation 
   \bes
u(t,x)=\int_\mathbb{R}p(t,x-y)u_0(y)dy+\lambda\int_0^t\int_\mathbb{R}p(t-s,x-y)
\sigma(u(s,y))w(dsdy),
  \lbl{4.7}\ees
where $p(t,x)$ is the transition density function of the symmetric $\alpha$-stable process $X_t$.

Before we state our main results, we recall some properties of the kernel function (transition density function)
$p(t,x)$.
\begin{prop}\lbl{p4.1} ({\rm\cite{S}}) The transition density $p(t,\cdot)$ of a strictly
$\alpha$-stable process satisfies

(i) $p(st,x)=t^{-1/\alpha}p\left(s,t^{-1/\alpha}x\right)$;

(ii) For $t$ large enough such that $p(t,0)\leq1$ and $a>2$, we have
 \bess
p(t,(x-y)/a)\geq p(t,x)p(t,y), \ \ \ {\rm for\ all }\ x\in\mathbb{R};
   \eess

(iii) $p(t,x)\asymp t^{-1/\alpha}\wedge\frac{t}{|x|^{1+\alpha}}$.
   \end{prop}
By using Proposition \ref{p4.1}, it is easy to verify that
   \bes
\int_\mathbb{R}p(t,x)p(s,x)dx=p(t+s,0).
   \lbl{4.8}\ees
In particular, $\|p(t,\cdot)\|_{L^2(\mathbb{R})}^2=p(2t,0)$.

Let
     \bess
\mathscr{E}_t(\lambda)=\sqrt{\mathbb{E}\left(\|u_t\|_{L^2(\mathbb{R})}^2\right)},
    \eess

   \begin{theo}\lbl{t4.2} Assume that {\rm(\ref{1.6})} holds and
   \bes
\int_\mathbb{R}\left(\int_\mathbb{R}p(t,x-y)u_0(y)dy\right)^2dx\leq \mu,
    \lbl{4.9}\ees
where $\mu$ is a positive constant,
then the noise excitation index of solution to {\rm(\ref{4.6})} with initial data $u_0(x)\geq0(\not\equiv0)$ is
$2\alpha/(\alpha-1)$.
  \end{theo}
\begin{remark}\lbl{r4.1} We remark the condition (\ref{4.9}) does indeed make sense.
Let us give an example. Taking $u_0(x)=\delta_{x_0}(x)$ for a fixed $x_0\in\mathbb{R}$, we have
     \bess
\int_\mathbb{R}\left(\int_\mathbb{R}p(t,x-y)u_0(y)dy\right)^2dx
=\int_\mathbb{R}p^2(t,x-x_0)dx
=p(2t,0)=:\mu<\infty. 
     \eess

Another example is that $u_0$ can be taken a function with compact support, such as the indicator  
function of a closed interval $[-l,l$ for arbitrarily fixed $l>0$, that is,  $u_0(x)=1_{[-l,l]}(x)$. Namely, 
$u_0(x)=1$ for $x\in[-l,l]$ and $u_0(x)=0$ for $x\not\in[-l,l]$. For a fixed $t>0$, it follows
from Proposition \ref{p4.1} that there is an $x_0>0$ such that
  \bess
p(t,x-y)\leq\frac{C}{|x-l|^{1+\alpha}} \ \ \ {\rm for}\ x>x_0,\,y\in[-l,l]; \\
p(t,x-y)\leq\frac{C}{|x+l|^{1+\alpha}} \ \ \ {\rm for}\ x<-x_0,\,y\in[-l,l],
   \eess
where $C>0$ is a constant. Direct calculations then show that
    \bess
\int_\mathbb{R}\left(\int_\mathbb{R}p(t,x-y)u_0(y)dy\right)^2dx
&=&\int_\mathbb{R}\left(\int_{-l}^lp(t,x-y)dy\right)^2dx\\
&\leq&\frac{C^2x_0}{t^{2/\alpha}}+2\int_{x_0}^\infty\frac{C^2}{|x-l|^{2+2\alpha}}dx
+2\int_{-\infty}^{-x_0}\frac{C^2}{|x+l|^{2+2\alpha}}dx\\
&=:&\mu(x_0,\alpha)<\infty.
  \eess
\end{remark}

\begin{lem}\lbl{l4.2}{\rm(\cite[Proposition 2.6]{MTL})} Let $T\leq\infty$ and $\beta>0$. Suppose 
that $f(t)$ is a nonnegative, locally integrable function satisfying
   \bes
f(t)\geq c_1+k\int_0^t(t-s)^{\beta-1}f(s)ds, \ \ \ \ {\rm for\ all}\ 0\leq t\leq T,
   \lbl{4.10}\ees
where $c_1$ is some positive constant. Then for any $t\in(0,T]$, we have the following
   \bess
\liminf_{k\rightarrow\infty}\frac{\log\log f(t)}{\log k}\geq\frac{1}{\beta}.
   \eess
\end{lem}
When the inequality (\ref{4.10}) is reversed with the second inequality in (\ref{4.3}), we have
\bess
\limsup_{k\rightarrow\infty}\frac{\log\log f(t)}{\log k}\leq\frac{1}{\beta}.
   \eess

{\bf Proof of Theorem \ref{t4.1}.} By using mild formulation and It\^{o} isometry,
we have the following
   \bes
\mathbb{E}|u(t,x)|^2=\Big|\int_\mathbb{R}p(t,x-y)u_0(y)dy\Big|^2
+\lambda^2\int_0^t\int_\mathbb{R}p^2(t-s,x-y)\mathbb{E}|\sigma(u(s,y))|^2dyds.
   \lbl{4.11}\ees
Integrating over $\mathbb{R}$ and by Fubini lemma, we get
\bes
\mathbb{E}\|u(t)\|_{L^2(\mathbb{R})}^2&=&\int_\mathbb{R}\Big|\int_\mathbb{R}
p(t,x-y)u_0(y)dy\Big|^2dx\nm\\ &&+\lambda^2\int_0^t\int_\mathbb{R}\mathbb{E}
|\sigma(u(s,y))|^2\left(\int_\mathbb{R}p^2(t-s,x-y)dx\right)dyds.
\lbl{4.12}\ees
By using (\ref{1.6}), (\ref{4.8}) and (\ref{4.9}), we have
\bes
\mathbb{E}\|u(t)\|_{L^2(\mathbb{R})}^2&\leq&\mu+\lambda^2\int_0^tp(2(t-s),0)\int_\mathbb{R}
\mathbb{E}|\sigma(u(s,y))|^2dyds\nm\\&\leq&\mu+\lambda^2L^2_\sigma\int_0^t\frac{C}{(t-s)^{1/\alpha}}
\mathbb{E}\|u(s)\|_{L^2(\mathbb{R})}^2ds.
   \lbl{4.13}\ees
By Lemma \ref{l4.1}, (\ref{4.13}) then implies that
   \bes
\limsup_{\lambda\rightarrow\infty}\frac{\log\log \mathbb{E}\|u(t)\|_{L^2(\mathbb{R})}^2}{\log \lambda}\leq\frac{2\alpha}{\alpha-1}.
   \lbl{4.14}\ees
That is, $\bar{\mathbbm{e}}(t)\leq 2\alpha/(\alpha-1)$.

Next, we prove $\underline{\mathbbm{e}}(t)\geq 2\alpha/(\alpha-1)$.
First, it follows from (\ref{4.11}) that
   \bess
\mathbb{E}|u(t,x)|^2\geq
\lambda^2\int_0^t\int_\mathbb{R}p^2(t-s,x-y)\mathbb{E}|\sigma(u(s,y))|^2dyds.
   \eess
Integrating over $\mathbb{R}$ with utilising Fubini lemma, we get
    \bes
\mathbb{E}\|u(t)\|_{L^2(\mathbb{R})}^2&\geq&\lambda^2\int_0^t\int_\mathbb{R}
\mathbb{E}|\sigma(u(s,y))|^2\left(\int_\mathbb{R}p^2(t-s,x-y)dx\right)dyds\nm\\
&\geq&\lambda^2l^2_\sigma\int_0^t\frac{C}{(t-s)^{1/\alpha}}\mathbb{E}\|u(s)\|_{L^2(\mathbb{R})}^2ds.
  \lbl{4.15}\ees
Again, by Lemma \ref{l4.1}, (\ref{4.15}) then implies that
   \bes
\liminf_{\lambda\rightarrow\infty}\frac{\log\log \mathbb{E}\|u(t)\|_{L^2(\mathbb{R})}^2}{\log \lambda}\geq\frac{2\alpha}{\alpha-1}.
   \lbl{4.16}\ees
That is, $\underline{\mathbbm{e}}(t)\geq 2\alpha/(\alpha-1)$.

Combining (\ref{4.14}) and (\ref{4.16}), we thus complete the proof of Theorem \ref{t4.1}. $\Box$

\begin{remark}\lbl{r4.2} When $\alpha=2$, $(\ref{4.16})$ was obtained in
Khoshnevisan-Kim {\rm\cite{KK1}}.
\end{remark}

\medskip


\begin{thebibliography}{99}\label{ref:ref}\addtolength{\itemsep}{-1.5ex}
\bibitem{BT} G. K. Batchelor and A. Townsend, {\em The nature of turbulent flow
at large wave numbers}, Proc. Royal Society A {\bf199} (1949) 238-255.

\bibitem{Cb} P-L. Chow, {\em Stochastic partial differential equations},
Chapman Hall/CRC Applied Mathematics and Nonlinear Science Series. Chapman Hall/CRC,
Boca Raton, FL, 2007. MR-2295103

\bibitem{C2009} P-L. Chow, {\em  Unbounded positive solutions of nonlinear parabolic It\^{o}
equations}, Communications on Stochastic Analysis {\bf 3} (2009) 211-222. MR-2588238

\bibitem{C2011} P-L. Chow, {\em  Explosive solutions of stochastic reaction-diffusion
equations in mean $L^p$-norm}, J. Differential Equations {\bf 250} (2011) 2567-2580. MR-2756076

\bibitem{R} R. C. Dalang, {\em  Extending the martingale measure stochastic integral with applications 
to spatially homogeneous s.p.d.e.'s}, Electron. J. Probab. {\bf 4} (1999) 1-29.  MR-1684157

\bibitem{DFV}  F. Delarue,  F. Flandoli and D. Vincenzi, {\em Noise prevents collapse of Vlasov-Poisson point charges},
Comm. Pure Appl. Math.  {\bf67}  (2014) 1700-1736. MR-3251910

\bibitem{DL} M. Dozzi and J.  L$\acute{o}$pez-Mimbela, {\em Finite-time blowup and existence of global positive
solutions of a semi-linear SPDE}, Stochastic Process. Appl. {\bf120}  (2010) 767-776. MR-2610325

\bibitem{E} H. W. Emmons, {\em The laminar-turbulent transition in a boundary layer Part 1},
J. Aeronaut. Sci {\bf18} (1951) 490-498. MR-0052931

\bibitem{FF} E. Fedrizzi and F. Flandoli, {\em Noise prevents singularities in linear transport equations}, J. Funct. Anal.  
{\bf264}  (2013), 1329-1354. MR-3017266

\bibitem{FGP}  F. Flandoli, M. Gubinelli and E. Priola, {\em Well-posedness of the transport equation by stochastic perturbation},
Invent. Math.  {\bf180}  (2010) 1-53.  MR-2593276

\bibitem{MM} M. Foondun and M. Joseph, {\em Remarks on non-linear noise excitability of some stochastic heat equations},
Stochastic Process. Appl. {\bf124} (2014) 3429-3440. MR-3231626

\bibitem{FK1} M. Foondun and D. Khoshnevisan, {\em  Intermittence and nonlinear parabolic stochastic partial differential 
equations}, Electron. J. Probab. {\bf 14} (2009) 548-568. MR-2480553 

\bibitem{GT} D. Gilbarg and N. S. Trudinger, {\em Elliptic partial differential equations of second order}, Classics in Mathematics, 
Springer-Verlag, Berlin, 2001. MR-1814364 

\bibitem{KK1} D. Khoshnevisan and K. Kim, {\em Non-linear noise excitation of intermittent stochastic PDEs and the topology of 
LCA group}, Ann. Probab. {\bf43} (2015) 1944-1991. MR-3353819

\bibitem{KK2} D. Khoshnevisan and K. Kim, {\em Non-linear noise excitation and intermittency under high disorder}, Proc. Amer. Math. 
Soc. {\bf143} (2015) 4073-4083. MR-3359595 

\bibitem{L} W. Liu, {\em Well-posedness of stochastic partial differential equations with Lyapunov condition}, J. Differential Equations 
{\bf 255} (2013) 572-592. MR-3053478

\bibitem{LiuRockner} W. Liu and M. Rockner, {\em Stochastic partial differential equations: an introduction}, Universitext. Springer, 
Cham, 2015. MR-3410409

\bibitem{MTL} W. Liu, K. Tian and M. Foondun, {\em On some properties of a class of fractional stochastic heat equations}, 
J. Theoret. Probab. {\bf 30} (2016) 1310-1333. MR-3736175 

\bibitem{LD} G. Y. Lv and J.  Duan, {\em Impacts of noise on a class of partial differential equations},
J. Differential Equations, {\bf258} (2015) 2196-2220. MR-3302534 

\bibitem{eP} E. Pardoux, {\em Stochastic partial differential equations and filtering of diffusion processes}, 
Stochastic, {\bf3} (1979) 127-167. MR-0553909

\bibitem{PZ} S. Peszat and J. Zabczyk, {\em Stochastic partial differential equations with
L\'evy noise (an evolution equation approach)}, Cambridge University Press, Cambridge, 2007. MR-2356959  

\bibitem{S} K. Sato, {\em L\'{e}vy processes and infinitely divisible distributions}, Cambridge Studies in Advanced Mathematics,
Vol. 68 Cambridge University Press, Cambridge, 1999. MR-1739520

\bibitem{Ta} T. Taniguchi, {\em The existence and uniqueness of energy solutions to local non-Lipschitz stochastic evolution equations}, 
J. Math. Anal. Appl. {\bf 360} (2009) 245-253. MR-2548380 

\bibitem{T} H. C. Tuckwell, {\em Stochastic process in the neurosciences}, CBMS-NSF
Regional Conference Series in Applied Mathematics, SIAM, Philadelphia, 1989. MR-1002192

\bibitem{Xie2016} B. Xie, {\em Some effects of the noise intensity upon non-linear stochastic heat equations on $[0,1]$}, Stochastic 
Process. Appl., {\bf 126} (2016) 1184-1205. MR-3461195

\end{thebibliography}
 \end{document}